\documentstyle{amsart}

\newcommand{\cN}{{\mathcal N}}
\newcommand{\fg}{{\mathfrak g}}

\newcommand{\BZ}{{\Bbb Z}}

\newcommand{\tr}{{\bf 1}}
\newcommand{\BC}{{\Bbb C}}
\newcommand{\cO}{{\mathcal O}}
\newcommand{\oC}{\overline{\cO}}

\title{Cohomology of subregular tilting modules for small quantum groups}
\thanks{The author is partially supported by the U.S. Civilian
 Research and Development Foundation under Award No. RM1-265.}
\address{Independent Moscow University, 11 Bolshoj Vlasjevskij per.,
 Moscow
121002 Russia}
\date{February, 1999}
\email{ostrik@@nw.math.msu.su}
\author{Viktor Ostrik}

\begin{document}
\maketitle
\section{Introduction}
Let $R$ be an irreducible root system with the Coxeter number $h$. 
Let $l>h$ be an odd integer
(we assume that $l$ is not divisible by 3 if $R$ is of
type $G_2$). Let $U$ be the quantum group of type 1 with divided powers
associated to these data, see \cite{Lu} (of type 1 means that the elements
$K_i^l$ are equal to 1). Let $u\subset U$ be the Frobenius kernel, see
{\em loc. cit.} Let $\tr$ be the trivial $U-$module. The cohomology
$H^{\bullet}(u,\tr)$ was computed by V.Ginzburg and S.Kumar in \cite{GK},
see also \cite{KV}. They proved
that the odd cohomology $H^{odd}(u,\tr)$ vanishes and the algebra of even
cohomology $H^{2\bullet}(u,\tr)$ is isomorphic to the algebra $\BC [\cN]$
of functions on the nilpotent cone $\cN \subset \fg$, where $\fg$ is the 
semisimple Lie algebra associated to $R$.
 Moreover, this is an isomorphism of graded
algebras with the grading on $\BC [\cN]$ corresponding to the
natural $\BC^*-$action on $\cN$ by dilatations. This isomorphism is
compatible with natural $G-$structures of both algebras where $G$ is
simply connected group associated to $R$.

Now let $s_a$ be the simple affine reflection lying in the affine Weyl group
associated
to $R,l$, see e.g. \cite{APW}. Let $\Theta_{s_a}$ be the corresponding
wall-crossing functor, see e.g. \cite{So}. Let $T=\Theta_{s_a}\tr$. It is
easy to see that cohomology $H^{\bullet}(u,T)$ has a natural algebra 
structure; namely for any simple $U-$module $L$ with highest weight lying
on the affine wall of the fundamental alcove we have $H^{\bullet}(u,T)=
Ext_u^{\bullet}(L,L)$. Since $T$ is a $U-$module the cohomology
$H^{\bullet}(u,T)$ has a natural structure of $G-$module.
 Let $\cO\subset \cN$ be the subregular nilpotent
orbit. The main result of this note is the following

{\bf Main Theorem.} {\em The odd cohomology $H^{odd}(u,T)$ vanishes.
The algebra $H^{2\bullet}(u,T)$ is isomorphic to the algebra
$\BC [ \oC ]$ of functions
on the closure of $\cO$. This is an 
isomorphism of graded algebras with the grading
on $\BC [\oC]$ corresponding to the action of $\BC^*$ by dilatations.
This isomorphism is compatible with natural $G-$structures of both
algebras.}

{\bf Remark.} One can prove the analogous theorem for the Frobenius kernel
$G_1$ of an almost simple algebraic
group $G$ over an algebraically closed field of characteristic $p>h$.

We remark that $\BC [\oC]=\BC [\cO]$ because of normality of $\oC$, see
\cite{Br, Th}. 

In \cite{He} W.H.Hesselink computed the structure of $\BC [\cN]$ as graded
$G$-module. It is easy to deduce the Hesselink Theorem from the
Ginzburg-Kumar Theorem (or rather from Andersen-Jantzen vanishing Theorem,
see \cite{AJ}). In the same way we are able to compute the
structure of $\BC[\oC]$ as graded $G-$module, see Corollary 3 below.

For any dominant weight $\lambda$ one defines the indecomposable tilting
module $T(\lambda)$ with highest weight $\lambda$, see e.g. \cite{An}.
Some time I believed that cohomology of any $T(\lambda)$ has parity
vanishing property. In fact this belief was main motivation for this work.
In the end of this note we give an example when cohomology of indecomposable
tilting module lives in both even and odd degrees.

\section{Proof of the Main Theorem}
Recall that $T$ has a 
unique trivial submodule $\tr$ and $T/\tr=H^0(s_a\cdot 0)$,
see e.g. \cite{An}. Let $\phi :T\to H^0(s_a\cdot 0)$ be the quotient map. 

{\bf Lemma 1.} {\em The map $\phi_*:H^{\bullet}(u,T)\to
H^{\bullet}(u,H^0(s_a\cdot 0))$ is zero.}

{\bf Proof.} The map $\phi_*$ is a map of $H^{2\bullet}(u,\tr)=\BC [\cN]$-
modules. It is known that the support of $H^{\bullet}(u,T)$ in $\cN$ is
equal to $\oC$, see \cite{O, J}. 

The cohomology
$H^{\bullet}(u,H^0(s_a\cdot 0))$ was computed by H.H.Andersen and J.C.Jantzen
in \cite{AJ}, 3.7. We reformulate their result as follows:

(a) Let $\pi :T^*(G/B)\to G/B$ be the cotangent bundle of the 
flag variety of the group
$G$. Let $s:T^*(G/B)\to \cN$ be the Springer
resolution. Let $L_{\theta}$ be the line bundle on $G/B$ corresponding
to the root $\theta$ dual to the highest coroot of $\fg$ (more directly
$\theta$ is the unique dominant short root). Then the even cohomology
$H^{ev}(u,H^0(s_a\cdot 0))$ vanishes; the odd cohomology is equal up to
shift to $s_*\pi^*L_{\theta}$ (if we consider the 
cohomology as a coherent sheaf on $\cN$).

In particular, if $\phi_*$ is nontrivial we obtain a section of the line bundle
$\pi^*L_{\theta}$ supported on $s^{-1}(\oC)$. Contradiction.

{\bf Remark.} In fact H.H.Andersen and J.C.Jantzen computed the cohomology
of induced modules over algebraic group over a field of characteristic
$p>0$. But their proof works in the quantum situation as well if we know some
vanishing result. This vanishing Theorem was proved in \cite{AJ} in types
$A,B,C,D,G$ or for strongly dominant weights. In our case the weight $\theta$
is not strongly dominant. A.Broer proved the desired vanishing in case of
characteristic 0 in \cite{Br}.
In recent work \cite{Th} all restrictions in Andersen-Jantzen
vanishing Theorem were removed. This should be used in the above-mentioned
generalization of our Main Theorem to characteristic $p$.

{\bf Corollary 1.} {\em The odd cohomology $H^{odd}(u,T)$ vanishes. For any
$i\ge 0$ we have an exact sequence
$$
0\to H^{2i-1}(u,H^0(s_a\cdot 0))\to H^{2i}(u,\tr)\to H^{2i}(u,T)\to 0.
$$
In particular, the natural map $H^{\bullet}(u,\tr)\to H^{\bullet}(u,T)$
is surjective.}

{\bf Proof} follows easily from consideration 
of the cohomology long exact sequence
associated with the short exact sequence
$$
0\to \tr \to T\to H^0(s_a\cdot 0)\to 0.
$$

{\bf Proof of the Main Theorem.} The surjectivity of the map 
$\BC [\cN]=H^{2\bullet}(u,\tr)\to H^{2\bullet}(u, T)$ implies that
there exists a surjection $\psi :H^{2\bullet}(u,T)\to \BC[\oC]$.
Let $L(\lambda)$ be the simple $G-$module with highest weight $\lambda$.
For any weight $\mu$ let $m_{\lambda}(\mu)$ be the multiplicity of the weight
$\mu$ in $L(\lambda).$ It is known that the multiplicity of $L(\lambda)$ in
$\BC [\oC ]$ is equal to $m_{\lambda}(0)-m_{\lambda}(\theta)$, see \cite{Br}
4.7. It is easy to deduce from Corollary 1 and (a) that the multiplicity of
$L(\lambda)$ in $H^{\bullet}(u,T)$ also equals
$m_{\lambda}(0)-m_{\lambda}(\theta)$ (we omit the proof since it is the same
as the proof of Corollary 3 below). Hence $\psi$ is an isomorphism.
The Theorem is proved.

Let $V=V(s_a\cdot 0)$ be the Weyl module with highest weight $s_a\cdot 0$.

{\bf Corollary 2.} {\em The cohomology $H^{\bullet}(u,V)$ is given by
$$
H^{2i}(u,V)=H^{2i}(u,T),
$$
$$
H^{2i+1}(u,V)=H^{2i}(u,\tr).
$$
}

{\bf Proof.} It is enough to consider the cohomology long exact sequence
associated with the short exact sequence
$$
0\to V\to T\to \tr \to 0
$$
and note that the map $H^{\bullet}(u,T)\to H^{\bullet}(u,\tr)$ is zero (this
can be proved in the same way as Lemma 1).

{\bf Remark.} One can easily compute the cohomology of the simple module
${\bf L}={\bf L}(s_a\cdot 0)$ with highest weight $s_a\cdot 0$ using the
short exact sequence
$$
0\to {\bf L}\to H^0(s_a\cdot 0)\to \tr \to 0.
$$
The answer is the following: $H^{2\bullet}(u,{\bf L})=0$ and for any
$i\ge 0$ we have short exact sequence
$$
0\to H^{2i}(u,\tr)\to H^{2i+1}(u,{\bf L})\to H^{2i+1}(u,H^0(s_a\cdot 0))\to 0.
$$

Let $R_+$ be the set of positive roots and let $W$ be the Weyl group.
For any $w\in W$ let $(-1)^w=det(w)$. Let $\rho$ be the halfsum of positive
roots. Let $w\cdot \lambda =w(\lambda +\rho)-\rho$.
For any dominant weight $\lambda$ let $d_n(\lambda)$ (resp. $t_n(\lambda)$) 
be the multiplicity of the simple module $L(\lambda)$ in the component of
degree $n$ of $\BC [\cN]$ (resp. $\BC [\oC]$). Let $p_n$ be the function on
the set $X$ of weights, given by
$$
\sum_{x\in X}\sum_{n\in \BZ}p_n(x)t^ne^x=\prod_{\alpha \in R_+}\frac{1}
{1-e^{\alpha}t}.
$$
This function is essentially the Kostant-Lusztig partition function.
Recall that Hesselink's Theorem (\cite{He}) states that $d_n(\lambda)=
\sum_{w\in W}(-1)^wp_n(w\cdot \lambda )$. Let $2k-1$ be the length of
reflection in $\theta$.

{\bf Corollary 3.} (cf. \cite{Br} 4.7) {\em We have
$$
t_n(\lambda)=\sum_{w\in W}(-1)^w(p_n(w\cdot \lambda )-p_{n-k}(w\cdot \lambda
-\theta)).
$$
}

{\bf Remark.} For types $A_l, B_l, C_l (l\ge 2), D_l (l\ge 3), G_2, F_4, E_6,
E_7,E_8$ the number $k$ equals to respectively $l, l, 2(l-1), 2l-3, 3, 8, 11,
17, 29$.

{\bf Proof.} Let $B$ be the Borel subgroup of $G$.
 Let $n$ be the nilpotent radical of the Borel subalgebra in
$\fg$. Let $S^{\bullet}(n^*)$ be the algebra of functions on $n$. By \cite{AJ,
Br} we have
$$
 H^{2i}(u,\tr)=Ind_B^G(S^i(n^*)), R^{>0}Ind_B^G(S^i(n^*))=0,
$$
$$
 H^{2i-1}(u,H^0(s_a\cdot 0))=Ind_B^G(S^{i-k}(n^*)\otimes \theta),
 R^{>0}Ind_B^G(S^{i-k}(n^*)\otimes \theta)=0.
$$
Now the Euler characteristic of $R^{\bullet}Ind_B^G(?)$ is given by the
Weyl character formula. The result follows.

{\bf Example.}
Here we present an example when cohomology (over
Frobenius kernel) of indecomposable
tilting module lives in both odd and even degrees.
Let $R$ be of type $A_2$. Let $s_1, s_2$
be the simple reflections in Weyl group, and let
$s_3$ be the affine reflection. Consider indecomposable tilting
module $T=T(s_3s_1s_2s_3\cdot 0)$. It has a filtration
with subquotients $H^0(s_3s_1s_2s_3\cdot 0), H^0(s_3s_1s_2\cdot 0),
H^0(s_3\cdot 0)$ and $H^0(0)$. Let $\omega_1$ and $\omega_2$
be the fundamental weights. We have $s_3s_1s_2s_3\cdot 0=(3l-3)\omega_2$.
 By Andersen-Jantzen Theorem
the cohomology of $H^0(s_3s_1s_2s_3\cdot 0)$ equals to
$Ind_B^G(3\omega_2\otimes S^{\bullet}(n^*))$ living in even degrees,
the cohomology of $H^0(s_3s_1s_2\cdot 0)$ or $H^0(s_3\cdot 0)$
equals to $Ind_B^G((\omega_1+\omega_2)\otimes S^{\bullet}(n^*))$
living in odd degrees, finally the cohomology of $H^0(0)$ equals
to $Ind_B^G(S^{\bullet}(n^*))$ living in even degrees.
By Kostant multiplicity formula we obtain that multiplicity of
$L(\lambda)$ in Euler characteristic of cohomology of $T$ equals to
$m_{\lambda}(3\omega_2)+m_{\lambda}(0)-2m_{\lambda}(\omega_1+\omega_2)$.
In particular, multiplicity of $L(0)$ equals to 1 and multiplicity
of $L(3\omega_1)$ equals to -1. This contradicts to the parity
vanishing.

{\large Acknowledgements.} I am grateful to M.Finkelberg for useful
conversations. Thanks also due to H.H.Andersen and J.Humphreys for valuable
suggestions. I would like to thank Aarhus University for
its hospitality while this note was being written.

\end{document}